\documentclass[12pt]{article}

\usepackage{mathrsfs}
\usepackage{amsfonts}
\usepackage{tikz}
\usepackage{amssymb,bm}
\usepackage{amsmath}
\usepackage{amsthm}
\usepackage{bbding}
\usetikzlibrary{matrix,shapes,snakes}

\allowdisplaybreaks
\newtheorem{thm}{Theorem}[section]

\newtheorem{lem}[thm]{Lemma}
\newtheorem{cor}[thm]{Corollary}
\newtheorem{prop}[thm]{Proposition}
\newtheorem{rmk}[thm]{Remark}

\def\proof{\vskip 1mm\noindent{\it Proof.}\quad}
\newcommand{\qedd}{\hspace*{\fill}$\Box$\medskip}

\def\det{\hbox{\rm{det\,}}}

\def\tr{\hbox{\rm{tr}}}
\def\Tr{\hbox{\rm{Tr}}}

\def\nor{{\rm{N}}}

\newcommand{\DF}[2]{{\displaystyle\frac{#1}{#2}}}

\begin{document}

\title{The compositional inverse of a class of linearized permutation
polynomials over $\mathbb{F}_{2^n}$,\\ $n$ odd}

\author{Baofeng Wu\thanks{Key Laboratory of Mathematics
Mechanization, AMSS, Chinese Academy of Sciences,
 Beijing 100190,  China. Email: wubaofeng@amss.ac.cn}}
 \date{}

\maketitle

\begin{abstract}
In this paper, the compositional inverses of a class of linearized
permutation polynomials of the form
$P(x)=x+x^2+\tr\left(\frac{x}{a}\right)$ over the finite field
$\mathbb{F}_{2^n}$ for an odd positive integer $n$ are explicitly
determined.\vskip .5em


\end{abstract}



\section{Introduction}\label{secintro}

Let $\mathbb{F}_{q}$ be the finite field with $q$ elements where $q$
is a prime or a prime power. A polynomial over $\mathbb{F}_{q}$ is
called a permutation polynomial (or a PP for short) if the map
\begin{eqnarray*}
  \mathbb{F}_{q} &\longrightarrow& \mathbb{F}_{q}\\
   c&\longmapsto&f(c)
\end{eqnarray*}
is bijective \cite{lidl}. It is clear that the set of all PPs over
$\mathbb{F}_{q}$ forms a group under the operation of composition of
polynomials and subsequent reduction modulo $(x^q-x)$, which is
isomorphic to $\mathcal {S}_q$, the symmetric group on $q$ letters.
Hence for every PP over $\mathbb{F}_{q}$, there exists a unique
polynomial $f^{-1}(x)$ over $\mathbb{F}_{q}$ such that
$f(f^{-1}(x))= f^{-1}(f(x))=x$ in the sense of reducing modulo
$(x^q-x)$, called the compositional inverse of $f(x)$. Generally
speaking, for a given PP $f(x)$ over $\mathbb{F}_{q}$, it is easy to
verify whether the other polynomial $g(x)$ is its compositional
inverse since we need only to check whether $f(g(x))\equiv x\mod
(x^q-x)$ or $g(f(x))\equiv x\mod (x^q-x)$ holds. But to derive the
compositional inverse of a PP explicitly is very challenging.
Indeed, it was posed by G.L. Mullen as a research problem that to
compute the coefficients of the compositional
 inverse  of a PP efficiently \cite[Problem
 10]{mullen1991}. Up to present, there are only a few classes of PPs
 over finite fields whose compositional inverses can be explicitly
 obtained (see, for example,  \cite{lidlDic,coulter,wu2}). For instance, a method
 to compute compositional inverses of linearized permutation
 polynomials was proposed by the authors in \cite{wu}.

Linearized polynomials are a special type of polynomials over finite
fields. A linearized polynomial over the finite field
$\mathbb{F}_{q^n}$ is of the form
\[L(x)=\sum_{i=0}^{n-1}a_ix^{q^i}.\]
It can induce a linear transformation of $\mathbb{F}_{q^n}$ over
$\mathbb{F}_{q}$ due to the well-known relation
$(x+y)^{q^t}=x^{q^t}+y^{q^t}$ for any $x,~y\in\mathbb{F}_{q^n}$ and
$0\leq t\leq n-1$. A criterion on linearized polynomials to be
permutation ones is easy to obtain. In fact, it was found by Dickson
as early as in 1897 that
$L(x)=\sum_{i=0}^{n-1}a_ix^{q^i}\in\mathbb{F}_{q^n}[x]$ is a PP if
and only if the matrix
\[D_L=\begin{pmatrix}
a_0&a_1&\dots&a_{n-1}\\
a_{n-1}^q&a_0^q&\dots&a_{n-2}^q\\
\vdots&\vdots&&\vdots\\
a_1^{q^{n-1}}&a_2^{q^{n-1}}&\dots&a_0^{q^{n-1}}
\end{pmatrix}\]
is non-singular \cite{lidl}. Obviously, when a linearized polynomial
is a PP, its compositional inverse  is also a linearized polynomial.

In \cite{wu}, $D_L$  is called the associate Dickson matrix of the
linearized polynomial $L(x)$, or vice versa, by the authors. We have
also found more relations between the properties of a linearized
polynomial and its associate Dickson matrix. For example, we have
related the compositional inverse of a linearized PP to the inverse
of its associate Dickson matrix, obtaining the following theorem:

\begin{thm}\label{invLPP}
Let $L(x)\in\mathbb{F}_{q^n}[x]$ be a linearized PP. Then
$D_{L^{-1}}=D_L^{-1}$. More precisely, if
$L(x)=\sum_{i=0}^{n-1}a_ix^{q^i}$ and $\bar{a}_i$ is the $(i,0)$-th
cofactor of $a_i$, $0\leq i\leq n-1$, then
\[L^{-1}(x)=\frac{1}{\det L}\sum_{i=0}^{n-1}\bar{a}_ix^{q^i},\]
where $\det
L=\sum_{i=0}^{n-1}a_{n-i}^{q^i}\bar{a}_i\in\mathbb{F}_{q}$
(subscripts reduced modulo $n$).
\end{thm}

According to Theorem \ref{invLPP}, the main task in computing the
compositional inverse of a linearized PP over $\mathbb{F}_{q^n}$ is
to compute the determinant of $n$ matrices over $\mathbb{F}_{q^n}$
with size $(n-1)\times(n-1)$, which can be efficiently realized from
an algorithmic perspective. However, to explicitly represent the
compositional inverses of certain linearized PPs by  closed form
formulas seems not so direct due to the difficulties in computing
determinants of computing matrices, even for some linearized PPs of
simple forms. In this paper, we devote to talking about such an
example.

Let $n$ be odd and $a\in\mathbb{F}_{2^n}$ be a non-zero element with
$\tr\left(\frac{1}{a}\right)=1$, where ``$\tr$" is the trace map
over $\mathbb{F}_{2^n}$, i.e. $\tr(x)=\sum_{i=0}^{n-1}x^{2^i}$ for
any $x\in\mathbb{F}_{2^n}$. It can be easily derived by
\cite[Theorem 7]{charpin} that
\[P(x)=x+x^2+\tr\left(\frac{x}{a}\right)\]
is a  linearized PP over $\mathbb{F}_{2^n}$. Linearized PPs of this
form can be related to the multiplication of the binary Knuth
pre-semifield in finite geometry \cite{biliotti}.

In the subsequent sections, we fix an odd $n$ and
$a\in\mathbb{F}_{2^n}$ with $\tr\left(\frac{1}{a}\right)=1$.  We
propose the explicit representation of coefficients of the
compositional inverse of $P(x)$ and explain how to obtain it
afterwards. It will be noted that $P^{-1}(x)$ is of a very
complicated form.


\section{Explicit representation of the compositional inverse of $P(x)$}\label{secinv}

To simplify the presentation of coefficients of $P^{-1}(x)$, we
introduce a symbol firstly.

For any $c\in\mathbb{F}_{2^n}$ and
$i_1,\ldots,i_t\in\mathbb{Z}/n\mathbb{Z}=\{0,1,\ldots,n-1\}$, we
denote
\[c\langle i_1,i_2,\ldots,i_t\rangle=\sum_{j=1}^{t}c^{2^{i_j}}.\]
Besides, for any subset $I$ of $\mathbb{Z}/n\mathbb{Z}$, we denote
\[c\langle I\rangle=\sum_{i\in I}c^{2^{i}}.\]
It is easy to see that the symbol ``$c\langle \cdot\rangle$"
satisfies the following properties:

\begin{prop}\label{symb1}
(1)
$c\langle0,1,\ldots,n-1\rangle=c\langle\mathbb{Z}\rangle=\tr(c)$;

\noindent(2) For any $\pi\in\mathcal{S}_t$,  $c\langle
i_{\pi(1)},i_{\pi(2)},\ldots,i_{\pi(t)}\rangle=c\langle
i_1,i_2,\ldots,i_t\rangle$;

\noindent(3) For any $ I,J\subseteq\mathbb{Z}/n\mathbb{Z}$,
$c\langle I\rangle+c\langle J\rangle=c\langle (I\cup
J)\backslash(I\cap J)\rangle$;

\noindent(4) For any $ k\geq 0$, $c\langle
i_1,i_2,\ldots,i_t\rangle^{2^k}=c\langle
i_1+k,i_2+k,\ldots,i_t+k\rangle$.\qedd
\end{prop}

By virtue of this newly defined symbol, we can propose the
compositional inverse of $P(x)$.

\begin{thm}\label{coeffinvPx}
$P^{-1}(x)=B_a(x)+\tr(x)$, where the coefficients of
$B_a(x)=\sum_{i=0}^{n-1}b_ix^{2^i}\in\mathbb{F}_{2^n}[x]$ are given
by
\begin{eqnarray*}
  b_0 &=& \frac{1}{a}\langle2,4,6,\ldots,n-1\rangle, \\
  b_i &=& \left\{\begin{array}{ll}
\frac{1}{a}\langle1,3,5,\ldots,i,\;i+1,i+3,\ldots,n-1\rangle& \text{if}~i~\text{is~odd},\\[.2cm]
\frac{1}{a}\langle1,3,5\ldots,i-1,\;i+2,i+4,\ldots,n-1\rangle&
\text{if}~i~\text{is~even},
  \end{array}
  \right.
\end{eqnarray*}
\noindent $1\leq i\leq n-1$го
\end{thm}
\proof Noting that
\begin{eqnarray*}
   B_a(1)&=&\sum_{i=0}^{n-1}b_i  \\
   &=&\DF{1}{a}\langle2,4,\ldots,n-1\rangle +\sum_{1\leq i\leq n-1\atop i~\text{odd}}
   \DF{1}{a}\langle1,3,\ldots,i,\;i+1,i+3,\ldots,n-1\rangle  \\
   && +\sum_{1\leq i\leq n-1\atop i~\text{even}}
   \DF{1}{a}\langle1,3\ldots,i-1,\;i+2,i+4,\ldots,n-1\rangle \\
   &=& \DF{1}{a}\langle2,4,\ldots,n-1\rangle +\sum_{j=1}^{\frac{n-1}{2}}\left(
   \DF{1}{a}\langle1,3,\ldots,2j-1,\;2j,2j+2,\ldots,n-1\rangle\right.\\
   &&\qquad\qquad\qquad\qquad~~\qquad+\left.\DF{1}{a}\langle1,3,\ldots,2j-1,\;2j+2,2j+4,\ldots,n-1\rangle\right) \\
   &=& \DF{1}{a}\langle2,4,\ldots,n-1\rangle +\sum_{j=1}^{\frac{n-1}{2}}\DF{1}{a}\langle2j\rangle \\
   &=&0,
\end{eqnarray*}
we get
\begin{eqnarray*}
  L^{-1}(L(x)) &=&B_a(x)+B_a(x^2)+B_a(1)\tr\left(\frac{x}{a}\right)+\tr(x+x^2)+\tr\left(\frac{x}{a}\right)  \\
   &=&B_a(x)+B_a(x^2)+\tr\left(\frac{x}{a}\right).
\end{eqnarray*}
Furthermore, by Proposition \ref{symb1} we have
\begin{eqnarray*}
  B_a(x)+B_a(x^2) &=&\DF{1}{a}\langle2,4,\ldots,n-1\rangle x+
  \sum_{1\leq i\leq n-1\atop i~\text{odd}}
   \DF{1}{a}\langle1,3,\ldots,i,\;i+1,i+3,\ldots,n-1\rangle x^{2^i}  \\
   &&+\sum_{1\leq i\leq n-1\atop i~\text{even}}
   \DF{1}{a}\langle1,3\ldots,i-1,\;i+2,i+4,\ldots,n-1\rangle x^{2^i} \\
   &&+\DF{1}{a}\langle2,4,\ldots,n-1\rangle x^2+
  \sum_{1\leq i\leq n-1\atop i~\text{odd}}
   \DF{1}{a}\langle1,3,\ldots,i,\;i+1,i+3,\ldots,n-1\rangle x^{2^{i+1}}  \\
   &&+\sum_{1\leq i\leq n-1\atop i~\text{even}}
   \DF{1}{a}\langle1,3\ldots,i-1,\;i+2,i+4,\ldots,n-1\rangle x^{2^{i+1}} \\
   &=&\left(\DF{1}{a}\langle2,4,\ldots,n-1\rangle+\DF{1}{a}\langle1,3,\ldots,n-2\rangle\right)x\\
   &&+\left(\DF{1}{a}\langle1,2,4,\ldots,n-1\rangle+\DF{1}{a}\langle2,4,\ldots,n-1\rangle\right)x^2\\
   &&+\sum_{3\leq i\leq n-1\atop
   i~\text{odd}}\left(\DF{1}{a}\langle1,3,\ldots,i,\;i+1,i+3,\ldots,n-1\rangle\right.\\
   &&\qquad\qquad~~+\left.\DF{1}{a}\langle1,3,\ldots,i-2,\;i+1,i+3,\ldots,n-1\rangle\right)x^{2^i}\\
   &&+\sum_{2\leq i\leq n-1\atop
   i~\text{even}}\left(\DF{1}{a}\langle1,3,\ldots,i-1,\;i+2,i+4,\ldots,n-1\rangle\right.\\
   &&\qquad\qquad~~+\left.\DF{1}{a}\langle1,3,\ldots,i-1,\;i,i+2,\ldots,n-1\rangle\right)x^{2^i}\\
   &=&\DF{1}{a}\langle1,2,3,4,\ldots,n-1\rangle x+\DF{1}{a}\langle1\rangle
   x^2+\sum_{2\leq i\leq n-1}\DF{1}{a}\langle i\rangle
   x^{2^i}\\
   &=&\left(\tr\left(\frac{1}{a}\right)+\frac{1}{a}\right)x+\sum_{1\leq i\leq
   n-1}\left(\frac{x}{a}\right)^{2^i}\\
   &=&x+\tr\left(\frac{x}{a}\right).
\end{eqnarray*}
Then $L^{-1}(L(x))=x$ follows. \qedd

\begin{rmk}
Assume $P^{-1}(x)=\sum_{i=0}^{n-1}p_ix^{2^i}$. Then it is easy to
derive from Theorem \ref{coeffinvPx} and the condition
$\tr\left(\frac{1}{a}\right)=1$ that
\begin{eqnarray*}
  p_0 &=& 1+\frac{1}{a}\langle2,4,\ldots,n-1\rangle \\
   &=&\frac{1}{a}\langle0,1,3,\ldots,n-2\rangle,
\end{eqnarray*}
and for $1\leq i\leq n-1$,
\begin{eqnarray*}
p_i &=& 1+\frac{1}{a}\langle1,3,\ldots,i,\;i+1,i+3,\ldots,n-1\rangle \\
   &=&\frac{1}{a}\langle0,2,4,\ldots,i-1,\;i+2,i+4,\ldots,n-2\rangle
\end{eqnarray*}
when $i$ is odd, and
\begin{eqnarray*}
  p_i &=& 1+\frac{1}{a}\langle1,3,\ldots,i-1,\;i+2,i+4,\ldots,n-1\rangle \\
   &=&\frac{1}{a}\langle0,2,4,\ldots,i,\;i+1,i+3,\ldots,n-2\rangle.
\end{eqnarray*}
when $i $ is even. Therefor, we can also represent $P^{-1}(x)$ by
$P^{-1}(x)=C_a(x)+\tr\left(\frac{x}{a}\right)$, where the
coefficients of
$C_a(x)=\sum_{i=0}^{n-1}c_ix^{2^i}\in\mathbb{F}_{2^n}[x]$ are given
by
\begin{eqnarray*}
  c_0 &=& \frac{1}{a}\langle1,3,5,\ldots,n-2\rangle, \\
  c_i &=& \left\{\begin{array}{ll}
1+\frac{1}{a}\langle1,3,5,\ldots,i-2,\;i+1,i+3,\ldots,n-1\rangle& \text{if}~i~\text{is~odd},\\[.2cm]
\DF{1}{a}\langle0,2,4\ldots,i-2,\;i+1,i+3,\ldots,n-2\rangle&
\text{if}~i~\text{is~even},
  \end{array}
  \right.
\end{eqnarray*}
\noindent $1\leq i\leq n-1$.
\end{rmk}

\begin{cor}\label{coeffinvP1x}
Let $P_1(x)=x+x^2+\tr(x)$. Then $P_1(x)$ is a linearized PP over
$\mathbb{F}_{2^n}$ with compositional inverse
\[P_1^{-1}(x)=\left\{\begin{aligned}
&\sum_{i=0}^{\frac{n-1}{2}}x^{2^{2i}}&~~\text{if}~~n\equiv1\mod4,\\
&\sum_{i=0}^{\frac{n-3}{2}}x^{2^{2i+1}}&~~\text{if}~~n\equiv 3\mod
4.
\end{aligned}
\right.\]
\end{cor}
\proof It is obvious that $P_1(x)$ is a linearized PP over
$\mathbb{F}_{2^n}$ since $\tr(1)=1$. By fixing $a=1$ in Theorem
\ref{coeffinvPx}, it is straightforward to obtain the coefficients
of $P_1^{-1}(x)$.\qedd

\begin{rmk}
(1) It is easy to see that the compositional inverse of $P_1(x)$ in
Corollary \ref{coeffinvP1x} can be written as
\[P_1^{-1}(x)=\frac{n+1}{2}\sum_{i=0}^{\frac{n-1}{2}}x^{2^{2i}}+\frac{n-1}{2}\sum_{i=0}^{\frac{n-3}{2}}x^{2^{2i+1}}.\]

\noindent (2) In fact, $P_1(x)$ is a linearized polynomial with
coefficients in $\mathbb{F}_{2}$, whose conventional associate
\cite{lidl} is $p_1(x)=\sum_{i=2}^{n-1}x^i$. By \cite[Lemma
3.59]{lidl}, it is easy to derive that the conventional associate of
$P_1^{-1}(x)$ is just the polynomial $\bar{p}(x)$ satisfying
$p(x)\bar{p}(x)\equiv1\mod(x^n+1)$. Hence from Corollary
\ref{coeffinvP1x}, we obtain
\[\bar{p}(x)=\left\{\begin{aligned}
&\sum_{i=0}^{\frac{n-1}{2}}x^{2i}&~~\text{if}~~n\equiv1\mod4,\\
&\sum_{i=0}^{\frac{n-3}{2}}x^{2i+1}&~~\text{if}~~n\equiv3\mod4.
\end{aligned}
\right.\]
\end{rmk}


\section{The method to obtain Theorem \ref{coeffinvPx}}\label{secmethod}

In this section, we describe in detail the process of deriving the
compositional inverse of $P(x)$. Note that the associate Dickson
matrix of $P(x)$ is
\[D_P=\begin{pmatrix}
1+\frac{1}{a}&1+\frac{1}{a^2}&\frac{1}{a^{2^2}}&\cdots&\frac{1}{a^{2^{n-2}}}&\frac{1}{a^{2^{n-1}}}\\
\frac{1}{a}&1+\frac{1}{a^2}&1+\frac{1}{a^{2^2}}&\cdots&\frac{1}{a^{2^{n-2}}}&\frac{1}{a^{2^{n-1}}}\\
\vdots&\vdots&\vdots&&\vdots&\vdots\\
\frac{1}{a}&\frac{1}{a^2}&\frac{1}{a^{2^2}}&\cdots&1+\frac{1}{a^{2^{n-2}}}&1+\frac{1}{a^{2^{n-1}}}\\
1+\frac{1}{a}&\frac{1}{a^2}&\frac{1}{a^{2^2}}&\cdots&\frac{1}{a^{2^{n-2}}}&1+\frac{1}{a^{2^{n-1}}}
\end{pmatrix},\]
which is of a relatively complicated form. To utilizing Theorem
\ref{invLPP} easily, we firstly simplify the problem by a linear
transformation of $P(x)$. That is, we reduce the  problem of
computing compositional inverse of $P(x)$ to that of computing
compositional inverse of $\tilde{P}(x)=ax+a^2x^2+\tr(x)$ by noting
that $\tilde{P}(x)=P(ax)$, which implies
\begin{equation}\label{eqn1}
P^{-1}(x)=a\tilde{P}^{-1}(x).
\end{equation}
In this place, we should deal with the Dickson matrix
\[D_{\tilde{P}}=\begin{pmatrix}
1+a&1+a^2&1&\cdots&1&1\\
1&1+a^2&1+a^{2^2}&\cdots&1&1\\
\vdots&\vdots&\vdots&&\vdots&\vdots\\
1&1&1&\cdots&1+a^{2^{n-2}}&1+a^{2^{n-1}}\\
1+a&1&1&\cdots&1&1+a^{2^{n-1}}
\end{pmatrix},\]
which seems much simpler than $D_P$.

Assume the $(i,0)$-th cofactor of $D_{\tilde{P}}$ is $\tilde{p}_i$,
$0\leq i\leq n-1$. Then
\[\det \tilde{P}=(1+a)\tilde{p}_0+\sum_{i=1}^{n-2}\tilde{p}_i+(1+a)\tilde{p}_{n-1}.\]
However, we can deduce $\det \tilde{P}=1$ before computing
$\tilde{p}_i$, $0\leq i\leq n-1$, since $\det \tilde{P}$ is a
non-zero element of $\mathbb{F}_2$ according to Theorem
\ref{coeffinvPx}. To explicitly represent $\tilde{p}_i$, $0\leq
i\leq n-1$, we give the following lemma on determinants of matrices
of a special form firstly.

\begin{lem}\label{detlemma1}
Let $a_1,\ldots,a_m,~b_1,\ldots,b_{m-1}\in R$ where $R$ is a
commutative ring with identity and $m\geq2$ is a positive integer.
Then
\[\det\begin{pmatrix}
a_1&b_1\\
&a_2&b_2\\
&&\ddots&\ddots\\
&&&a_{m-1}&b_{m-1}\\
1&1&\cdots&1&a_m
\end{pmatrix}=(-1)^{m+1}\prod_{j=1}^{m-1}b_j+\sum_{i=1}^{m-2}(-1)^{m+i+1}\prod_{j=1}^{i}\prod_{k=i+1}^{m-1}a_jb_k
+\prod_{j=1}^{m}a_j.\]
\end{lem}
\proof The determinant can be directly computed using Laplace
expansion along the last row of the matrix.\qedd

To simplify the computations of $\tilde{p}_i$, $0\leq i\leq n-1$, we
introduce another symbol and propose several obvious properties of
it.

For any $c\in\mathbb{F}_{2^n}^*$,
$i_1,\ldots,i_t\in\mathbb{Z}/n\mathbb{Z}=\{0,1,\ldots,n-1\}$ and
$I\subseteq\mathbb{Z}/n\mathbb{Z}$, we denote
\[c[ i_1,i_2,\ldots,i_t]=\prod_{j=1}^{t}c^{2^{i_j}}=c^{2^{i_1}+\cdots+2^{i_t}}\]
and
\[c[I]=\prod_{i\in I}c^{2^{i}}.\]
\begin{prop}\label{symb2}
(1) $c[0,1,\ldots,n-1]=\nor(c)=1$, where ``$\nor$" is the norm map
over $\mathbb{F}_{2^n}$;

\noindent(2) For any $\pi\in\mathcal{S}_t$, $c[
i_{\pi(1)},i_{\pi(2)},\ldots,i_{\pi(t)}]=c[ i_1,i_2,\ldots,i_t]$;

\noindent(3) For any $ I,~J\subseteq\mathbb{Z}/n\mathbb{Z}$, $c[
I]c[ J]=c[I\cup J]c[I\cap J]=c[ (I\cup J)\backslash(I\cap J)]c[I\cap
J]^2$;

\noindent(4) For any $ k\geq 0$, $c[ i_1,i_2,\ldots,i_t]^{2^k}=c[
i_1+k,i_2+k,\ldots,i_t+k]$;

\noindent(5) For any $i\subseteq\mathbb{Z}/n\mathbb{Z}$,
$c[i]=c\langle i\rangle$.\qedd
\end{prop}

Based on these preparations, $\tilde{p}_i$ are computed in the
following lemmas for $0\leq i\leq n-1$.

\begin{lem}\label{detd0}
\[\tilde{p}_0=\frac{1}{a}\left(1+\frac{1}{a}\langle2,4,\ldots,n-1\rangle\right).\]
\end{lem}
\proof
\begin{eqnarray*}
  \tilde{p}_0 &=& \det\begin{pmatrix}
1+a[1]&1+a[2]&1&\cdots&1&1\\
1&1+a[2]&1+a[3]&\cdots&1&1\\
\vdots&\vdots&\vdots&&\vdots&\vdots\\
1&1&1&\cdots&1+a[n-2]&1+a[n-1]\\
1&1&1&\cdots&1&1+a[n-1]
  \end{pmatrix} \\
   &=&\det\begin{pmatrix}
   a[1]&&a[3]\\
   &a[2]&&a[4]\\
   &&\ddots&&\ddots\\
   &&&a[n-3]&&a[n-1]\\
   &&&&a[n-2]\\
   1&1&\cdots&1&1&1+a[n-1]
   \end{pmatrix}.
\end{eqnarray*}
By Laplace expansion we get
\[\tilde{p}_0 =a[1]a[3]\cdots a[n-4]a[n-2]\cdot\det\begin{pmatrix}
a[2]&a[4]\\
&a[4]&a[6]\\
&&\ddots&\ddots\\
&&&a[n-3]&a[n-1]\\
1&1&\cdots&1&1+a[n-1]
\end{pmatrix}.\]
Then from Lemma \ref{detlemma1} and Proposition \ref{symb2} we
obtain
\begin{eqnarray*}
   \tilde{p}_0&=& a[1,3,\ldots,n-2]\cdot \Bigg[a[4,6,\ldots,n-1]+
\sum_{j=1}^{\frac{n-5}{2}}a[2,4,\ldots,2j]a[2j+4,2j+6,\ldots,n-1]\\
&&\qquad\qquad\qquad\qquad~~+a[2,4,\ldots,n-3](1+a[n-1])
   \Bigg]\\
   &=& \DF{\nor(a)}{a[0,2]}+ \sum_{j=1}^{\frac{n-5}{2}}\DF{\nor(a)}{a[0,2j+2]}+\DF{\nor(a)}{a[0,n-1]}(1+a[n-1])\\
   &=&\DF{1}{a}\left(\DF{1}{a[2]}+\sum_{j=1}^{\frac{n-5}{2}}\DF{1}{a[2j+2]}+\DF{1}{a[n-1]}+1\right)\\
   &=&\frac{1}{a}\left(1+\frac{1}{a}\langle2,4,\ldots,n-1\rangle\right).
\end{eqnarray*}\qedd

\begin{lem}\label{detd1}
\[\tilde{p}_1=\frac{1}{a}\left(1+\frac{1}{a}\langle1,2,4,\ldots,n-1\rangle\right).\]
\end{lem}
\proof
\begin{eqnarray*}
  \tilde{p}_1 &=& \det\begin{pmatrix}
1+a[1]&1&1&\cdots&1&1\\
1&1+a[2]&1+a[3]&\cdots&1&1\\
\vdots&\vdots&\vdots&&\vdots&\vdots\\
1&1&1&\cdots&1+a[n-2]&1+a[n-1]\\
1&1&1&\cdots&1&1+a[n-1]
  \end{pmatrix} \\
   &=&\det\begin{pmatrix}
   a[1]&a[2]&a[3]\\
   &a[2]&&a[4]\\
   &&a[3]&&a[5]\\
   &&&\ddots&&\ddots\\
  & &&&a[n-3]&&a[n-1]\\
  & &&&&a[n-2]\\
  1& 1&1&\cdots&1&1&1+a[n-1]
   \end{pmatrix}.
\end{eqnarray*}
By Laplace expansion we get
\[\tilde{p}_1 =a[3]a[5]\cdots a[n-2]\cdot\det\begin{pmatrix}
a[1]&a[2]\\
&a[2]&a[4]\\
&&\ddots&\ddots\\
&&&a[n-3]&a[n-1]\\
1&1&\cdots&1&1+a[n-1]
\end{pmatrix}.\]
Then from Lemma \ref{detlemma1} and Proposition \ref{symb2} we
obtain
\begin{eqnarray*}
   \tilde{p}_1&=& a[3,5,\ldots,n-2]\cdot
   \Bigg[a[2,4,\ldots,n-1]+a[1,4,6,\ldots,n-1]\\
&&+\sum_{j=2}^{\frac{n-3}{2}}a[1,2,4,\ldots,2j-2]a[2j+2,2j+4,\ldots,n-1]\\
&& +a[1,2,4,\ldots,n-3](1+a[n-1])
   \Bigg]\\
   &=& \DF{\nor(a)}{a[0,1]}+\DF{\nor(a)}{a[0,2]}+ \sum_{j=2}^{\frac{n-3}{2}}\DF{\nor(a)}{a[0,2j]}+\DF{\nor(a)}{a[0,n-1]}(1+a[n-1])\\
   &=&\DF{1}{a}\left(\DF{1}{a[1]}+\DF{1}{a[2]}+\sum_{j=2}^{\frac{n-3}{2}}\DF{1}{a[2j]}+\DF{1}{a[n-1]}+1\right)\\
   &=&\frac{1}{a}\left(1+\DF{1}{a}\langle1,2,4,\ldots,n-1\rangle\right).
\end{eqnarray*}\qedd

\begin{lem}\label{detdi}
For $2\leq i\leq n-3$,
\[\tilde{p}_i=\frac{1}{a}\left(1+\DF{1}{a}\langle1,3,\ldots,i-1,i+2,i+4,\ldots,n-1\rangle\right)\]
when $i$ is even, and
\[\tilde{p}_i=\frac{1}{a}\left(1+\DF{1}{a}\langle1,3,\ldots,i,i+1,i+3,\ldots,n-1\rangle\right)\]
when $i$ is odd.
\end{lem}
\proof
\begin{eqnarray*}
  \tilde{p}_i &=& \det\tiny \begin{pmatrix}
1+a[1]&1&\cdots&1&1&1&1&\cdots&1&1\\
1+a[1]&1+a[2]&\cdots&1&1&1&1&\cdots&1&1\\
\vdots&\vdots&&\vdots&\vdots&\vdots&\vdots&&\vdots&\vdots\\
1&1&\cdots&1+a[i-1]&1+a[i]&1&1&\cdots&1&1\\
1&1&\cdots&1&1&1+a[i+1]&1+a[i+2]&\cdots&1&1\\
\vdots&\vdots&&\vdots&\vdots&\vdots&\vdots&&\vdots&\vdots\\
1&1&\cdots&1&1&1&1&\cdots&1+a[n-2]&1+a[n-1]\\
1&1&\cdots&1&1&1&1&\cdots&1&1+a[n-1]
  \end{pmatrix} \\
   &=&\normalsize\det\scriptsize \begin{pmatrix}
   &a[2]\\
   a[1]&&a[3]\\
   &\ddots&&\ddots\\
   &&a[i-2]&&a[i]\\
   &&&a[i-1]&a[i]&a[i+1]&a[i+2]\\
   &&&&&a[i+1]&&a[i+3]\\
   &&&&&&\ddots&&\ddots\\
   &&&&&&&a[n-3]&&a[n-1]\\
   &&&&&&&&a[n-2]\\
   1&\cdots&1&1&1&1&\cdots&1&1&1+a[n-1]
   \end{pmatrix}.
\end{eqnarray*}\normalsize

\noindent(1) When $i$ is even, we can get by Laplace expansion that
\begin{eqnarray*}
  \tilde{p}_i &=& a[2,4,\ldots,i,i+1,i+3,\ldots,n-2]\cdot \\
   &&\det\small\begin{pmatrix}
a[1]&a[3]\\
&a[3]&a[5]\\
&&\ddots&\ddots\\
&&&a[i-3]&a[i-1]\\
&&&&a[i-1]&a[i+2]\\
&&&&&a[i+2]&a[i+4]\\
&&&&&&\ddots&\ddots\\
&&&&&&&a[n-3]&a[n-1]\\
1&1&\cdots&1&1&1&\cdots&1&1+a[n-1]
\end{pmatrix}.
\end{eqnarray*}\normalsize
Then from Lemma \ref{detlemma1} and Proposition \ref{symb2} we
obtain
\begin{eqnarray*}
   \tilde{p}_i&=& a[2,4,\ldots,i,i+1,i+3,\ldots,n-2]\cdot
   \Bigg[a[3,5,\ldots,i-1,i+2,i+4,\ldots,n-1]\\
&&+\sum_{j=2}^{\frac{n-4}{2}}a[1,3,\ldots,2j-1,2j+3,\ldots,i-1,i+2,i+4,\ldots,n-1]\\
&&+a[1,3,\ldots,i-3,i+2,i+4,\ldots,n-1]+a[1,3,\ldots,i-1,i+4,i+6,\ldots,n-1]\\
&&+\sum_{j=2}^{\frac{n-i-3}{2}}a[1,3,\ldots,i-1,i+2,i+4,\ldots,i+2j-2,i+2j+2\ldots,n-1]\\
 && +a[1,3,\ldots,i-1,i+2,i+4,\ldots,n-3](1+a[n-1])
   \Bigg]\\
   &=& \DF{\nor(a)}{a[0,1]}+
   \sum_{j=2}^{\frac{n-4}{2}}\DF{\nor(a)}{a[0,2j+1]}+\DF{\nor(a)}{a[0,i-1]}+\DF{\nor(a)}{a[0,i+2]}\\
   &&+\sum_{j=2}^{\frac{n-i-3}{2}}\DF{\nor(a)}{a[0,2j]}+\DF{\nor(a)}{a[0,n-1]}(1+a[n-1])\\
   &=&\frac{1}{a}\left(1+\DF{1}{a}\langle1,3,\ldots,i-1,i+2,i+4,\ldots,n-1\rangle\right);
\end{eqnarray*}

\noindent(2) When $i$ is odd, we can get by Laplace expansion that
\begin{eqnarray*}
  \tilde{p}_i &=& a[2,4,\ldots,i-1,i+2,i+4,\ldots,n-2]\cdot \\
   &&\det\small\begin{pmatrix}
a[1]&a[3]\\
&a[3]&a[5]\\
&&\ddots&\ddots\\
&&&a[i-2]&a[i]\\
&&&&a[i]&a[i+1]\\
&&&&&a[i+1]&a[i+3]\\
&&&&&&\ddots&\ddots\\
&&&&&&&a[n-3]&a[n-1]\\
1&1&\cdots&1&1&1&\cdots&1&1+a[n-1]
\end{pmatrix}.
\end{eqnarray*}\normalsize
Then from Lemma \ref{detlemma1} and Proposition \ref{symb2} we
obtain
\begin{eqnarray*}
   \tilde{p}_i&=& a[2,4,\ldots,i-1,i+2,i+4,\ldots,n-2]\cdot
   \Bigg[a[3,5,\ldots,i,i+1,i+3,\ldots,n-1]\\
&&+\sum_{j=2}^{\frac{n-3}{2}}a[1,3,\ldots,2j-1,2j+3,\ldots,i,i+1,i+3,\ldots,n-1]\\
&&+a[1,3,\ldots,i-2,i+1,i+3,\ldots,n-1]+a[1,3,\ldots,i,i+3,i+5,\ldots,n-1]\\
&&+\sum_{j=2}^{\frac{n-i-2}{2}}a[1,3,\ldots,i,i+1,i+3,\ldots,i+2j-3,i+2j+1\ldots,n-1]\\
 && +a[1,3,\ldots,i,i+1,i+3,\ldots,n-3](1+a[n-1])
   \Bigg]\\
   &=& \DF{\nor(a)}{a[0,1]}+
   \sum_{j=2}^{\frac{n-3}{2}}\DF{\nor(a)}{a[0,2j+1]}+\DF{\nor(a)}{a[0,i]}+\DF{\nor(a)}{a[0,i+1]}\\
   &&+   \sum_{j=2}^{\frac{n-i-2}{2}}\DF{\nor(a)}{a[0,i+2j-1]}+\DF{\nor(a)}{a[0,n-1]}(1+a[n-1])\\
   &=&\frac{1}{a}\left(1+\DF{1}{a}\langle1,3,\ldots,i,i+1,i+3,\ldots,n-1\rangle\right).
\end{eqnarray*}\qedd

\begin{lem}\label{detdn-2}
\[\tilde{p}_{n-2}=\frac{1}{a}\left(1+\DF{1}{a}\langle1,3,\ldots,n-2,n-1\rangle\right).\]
\end{lem}
\proof
\begin{eqnarray*}
  \tilde{p}_{n-2} &=& \det\begin{pmatrix}
1+a[1]&1&\cdots&1&1&1\\
1+a[1]&1+a[2]&\cdots&1&1&1\\
\vdots&\vdots&&\vdots&\vdots&\vdots\\
1&1&\cdots&1+a[n-3]&1+a[n-2]&1\\
1&1&\cdots&1&1&1+a[n-1]
  \end{pmatrix} \\
   &=&\det\begin{pmatrix}
   &a[2]\\
   a[1]&&a[3]\\
   &a[2]&&a[4]\\
   &&\ddots&&\ddots\\
   &&&a[n-4]&&a[n-2]\\
   &&&&a[n-3]&a[n-2]&a[n-1]\\
   1&1&\cdots&1&1&1&1+a[n-1]
   \end{pmatrix}.
\end{eqnarray*}
By Laplace expansion we get
\[\tilde{p}_{n-2} =a[2,4,\ldots,n-3]\cdot\det\begin{pmatrix}
a[1]&a[3]\\
&a[3]&a[5]\\
&&\ddots&\ddots\\
&&&a[n-4]&a[n-2]\\
&&&&a[n-2]&a[n-1]\\
1&1&\cdots&1&1&1+a[n-1]
\end{pmatrix}.\]
Then from Lemma \ref{detlemma1} and Proposition \ref{symb2} we
obtain
\begin{eqnarray*}
   \tilde{p}_{n-2}&=& a[2,4,\ldots,n-3]\cdot \Bigg[a[3,5,\ldots,n-2,n-1]\\
   &&+\sum_{j=1}^{\frac{n-5}{2}}a[1,3,\ldots,2j-1,2j+3,\ldots,n-2,n-1]\\
&&+a[1,3,\ldots,n-4,n-1]+a[1,3,\ldots,n-2](1+a[n-1])
   \Bigg]\\
   &=& \DF{\nor(a)}{a[0,1]}+ \sum_{j=1}^{\frac{n-5}{2}}\DF{\nor(a)}{a[0,2j+1]}+\DF{\nor(a)}{a[0,n-2]}+\DF{\nor(a)}{a[0,n-1]}(1+a[n-1])\\
   &=&\frac{1}{a}\left(1+\DF{1}{a}\langle1,3,\ldots,n-2,n-1\rangle\right).
\end{eqnarray*}\qedd

\begin{lem}\label{detdn-1}
\[\tilde{p}_{n-1}=\frac{1}{a}\left(1+\DF{1}{a}\langle1,3,\ldots,n-2\rangle\right).\]
\end{lem}
\proof
\begin{eqnarray*}
  \tilde{p}_{n-1} &=& \det\begin{pmatrix}
1+a[1]&1&\cdots&1&1\\
1+a[1]&1+a[2]&\cdots&1&1\\
\vdots&\vdots&&\vdots&\vdots\\
1&1&\cdots&1+a[n-2]&1+a[n-1]\\
  \end{pmatrix} \\
   &=&\det\begin{pmatrix}
   &a[2]\\
   a[1]&&a[3]\\
   &a[2]&&a[4]\\
   &&\ddots&&\ddots\\
   &&&a[n-3]&&a[n-1]\\
   1&1&\cdots&1&1+a[n-2]&1+a[n-1]
   \end{pmatrix}.
\end{eqnarray*}
By Laplace expansion we get
\[\tilde{p}_{n-1} =a[2,4,\ldots,n-1]\cdot\det\begin{pmatrix}
a[1]&a[3]\\
&a[3]&a[5]\\
&&\ddots&\ddots\\
&&&a[n-4]&a[n-2]\\
1&1&\cdots&1&1+a[n-2]
\end{pmatrix}.\]
Then from Lemma \ref{detlemma1} and Proposition \ref{symb2} we
obtain
\begin{eqnarray*}
   \tilde{p}_{n-1}&=& a[2,4,\ldots,n-1]\cdot \Bigg[a[3,5,\ldots,n-2]
   +\sum_{j=1}^{\frac{n-5}{2}}a[1,3,\ldots,2j-1,2j+3,\ldots,n-2]\\
   &&\qquad\qquad\qquad\qquad+a[1,3,\ldots,n-4](1+a[n-2])
   \Bigg]\\
   &=& \DF{\nor(a)}{a[0,1]}+ \sum_{j=1}^{\frac{n-5}{2}}\DF{\nor(a)}{a[0,2j+1]}+\DF{\nor(a)}{a[0,n-2]}(1+a[n-2])\\
   &=&\frac{1}{a}\left(1+\DF{1}{a}\langle1,3,\ldots,n-2\rangle\right).
\end{eqnarray*}\qedd

\vspace{\baselineskip}

From Lemma \ref{detd0}--\ref{detdn-1}, we can get the explicit
representation of $\tilde{P}^{-1}(x)$. Finally, Theorem
\ref{coeffinvPx} follows from \eqref{eqn1}.



\begin{thebibliography}{5}


\bibitem{biliotti}
M. Biliotti, V. Jha, N. Johnson, Foundations of translation planes,
Pure and Applied Mathematics, vol. 253,  Marcel Dekker, New York,
2001.

\bibitem{charpin}
P. Charpin, G. Kyureghyan, When does $G(x)+\gamma \Tr(H(x))$
permutate $\mathbb{F}_{p^n}$?, Finite Fields Appl. 15 (2009)
615--632.

\bibitem{coulter}
R.S. Coulter, M. Henderson, The compositional inverse of a class of
permutation polynomials over a finite field, Bull. Austral. Math.
Soc. 65 (2002) 521--526.


\bibitem{lidl}
R. Lidl, H. Niederreiter, Finite fields, second edn., Encyclopedia
Math. Appl., vol. 20, Cambridge University Press, Cambridge, 1997.

\bibitem{lidlDic}
R. Lidl, G.L. Mullen, G. Turnwald, Dickson polynomials, Pitman
monographs and surveys in pure and applied mathematics, vol. 65,
Longman Scientific \& Technical, Essex,  1993.


\bibitem{mullen1991}
G.L. Mullen, Permutation polynomials over finite fields, in: Finite
Fields, Coding Theory, and Advances in Communication and Computing,
Las Vegas, NY, 1991, pp. 131--151.




\bibitem{wu}
B. Wu, Z. Liu, Linearized polynomials over finite fields revisited,
Finite Fields Appl. 22 (2013) 79--100.

\bibitem{wu2}
B. Wu, Z. Liu, The compositional inverse of a class of bilinear
permutation polynomials over finite fields of characteristic 2,
Finite Fields Appl. (2013), doi: 10.1016/j.ffa.2013.05.003

\end{thebibliography}
\end{document}